\documentclass{article}
\usepackage{etoolbox}

\newtoggle{ARXIV}
\toggletrue{ARXIV}

\usepackage[usenames,dvipsnames]{xcolor}
\usepackage{graphicx, pgfplots, subfigure, standalone, pgf, tikz}
\usepackage[justification=centering]{caption}

\usepackage{comment, url,float, setspace, xspace, pifont}
\pgfplotsset{compat=1.18}

\usepackage[inline]{enumitem}
\usepackage{array, bm}
\setlist[enumerate]{label=(\roman*)}

\usepackage[ruled,vlined]{algorithm2e}
\usepackage[colorlinks=true,breaklinks=true,bookmarks=false,urlcolor=blue, citecolor=blue,linkcolor=blue,bookmarksopen=false,draft=false]{hyperref}
\iftoggle{ARXIV}{
\usepackage{amsfonts,amsthm,amsmath}
}{}
\usepackage[nameinlink]{cleveref}

\usepackage{csvsimple, booktabs, longtable, multirow, multicol, adjustbox}	

\usepackage[xindy, acronyms, nowarn]{glossaries}   
\makeglossaries
\newcommand{\eg}{{\it e.g.}}
\newcommand{\ie}{{\it i.e.}}
\newcommand{\etal}{{\it et al.}}

\definecolor{lime}{HTML}{A6CE39}
\DeclareRobustCommand{\orcidicon}{
	\begin{tikzpicture} \draw[lime, fill=lime] (0,0) circle [radius=0.16] node[white] { {\fontfamily{qag}\selectfont \tiny ID} };
	\draw[white, fill=white] (-0.0625,0.095) circle [radius=0.007];
	\end{tikzpicture} \hspace{-2mm}
}

\foreach \x in {A, ..., Z}{\expandafter\xdef\csname orcid\x\endcsname{\noexpand\href{https://orcid.org/\csname orcidauthor\x\endcsname}{\noexpand\orcidicon} } }
\newcommand{\maxi}{\text{maximize}}
\newcommand{\xofi}{x^i}
\newcommand{\xminusi}{x^{-i}}
\newcommand{\X}[1]{ \mathcal{X}^{ #1 } }

\newcommand{\NPhard}{\mbox{$\mathcal{NP}$-hard}\xspace}
\newcommand{\SigmaTwoP}{\mbox{$\Sigma^p_2$-complete}\xspace}

\newcommand{\conv}{\operatorname{conv}}

\newcommand{\cl}{\operatorname{cl}}

\newcommand{\cmark}{\ding{51}} 
\newcommand{\xmark}{\ding{55}} 

\crefname{line}{Line}{Lines}
\crefname{lemma}{Lemma}{Lemmata}
\crefname{theorem}{Theorem}{Theorems}
\crefname{proposition}{Proposition}{Propositions}
\crefname{algorithm}{Algorithm}{Algorithms}
\crefname{equation}{Equation}{Equations}
\crefname{definition}{Definition}{Definition}
\crefname{claim}{Claim}{Claim}
\crefname{corollary}{Corollary}{Corollaries}
\crefname{remark}{Remark}{Remarks}
\crefname{example}{Example}{Examples}
\crefname{figure}{Figure}{Figures}
\crefname{section}{Section}{Sections}
\crefname{table}{Table}{Tables}

\newacronym{MIP}{MIP}{Mixed-Integer Programming}
\newacronym{AGT}{AGT}{Algorithmic Game Theory}
\newacronym{POA}{PoA}{Price of Anarchy}
\newacronym{POS}{PoS}{Price of Security}
\newacronym{MNE}{MNE}{Mixed Nash Equilibrium}
\newglossaryentry{IPG}
{
  name={IPG},
  description={Integer Programming Game},
  first={\glsentrydesc{IPG} (\glsentrytext{IPG})},
  plural={IPGs},
  descriptionplural={Integer Programming Games},
  firstplural={Integer Programming Games (\glsentryplural{IPG})}
}

\newacronym{SGM}{SGM}{Sample Generation Method}
\newacronym{eSGM}{eSGM}{Exhaustive Sample Generation Method}
\newacronym{CNP}{CnP}{Cut-and-Play}
\newacronym{ZEROR}{ZR}{Zero Regrets}
\newacronym{BNP}{BnP}{Branch-and-Prune}
\newacronym{BM}{BM}{Branching Method}
\newacronym{GBM}{gBM}{Generalized Branching Method}
\newacronym{CNG}{CNG}{Critical Node Game}
\newacronym{MCNP}{MCNP}{Multilevel Critical Node Problem}

\newcommand{\defender}{f^{\mathtt{d}}}
\newcommand{\attacker}{f^{\mathtt{a}}}
\newcommand{\pprotect}{p^{\mathtt{d}}}
\newcommand{\pattack}{p^{\mathtt{a}}}

\iftoggle{ARXIV}{
    \usepackage[sort,numbers]{natbib}
    \usepackage[margin=1.2in]{geometry}
    \usepackage[justification=centering]{caption}
    \usepackage[title]{appendix}
    \usepackage{charter}

    \newtheorem{definition}{Definition}
    
    \newtheorem{example}{Example}
    
}

\begin{document}

\newcommand{\orcidauthorA}{0000-0002-2344-0960}
\newcommand{\orcidauthorB}{0000-0002-5447-6245}
\newcommand{\orcidauthorC}{0000-0001-9269-633X}
\newcommand{\orcidauthorD}{0000-0002-4662-3241}

\newcommand{\theTitle}{Integer Programming Games: A Gentle Computational Overview}
\newcommand{\theAbstract}{%
In this tutorial, we present a computational overview on computing Nash equilibria in Integer Programming Games (\emph{IPG}s), \ie, how to compute solutions for a class of non-cooperative and nonconvex games where each player solves a mixed-integer optimization problem. \emph{IPGs} are a broad class of games extending the modeling power of mixed-integer optimization to multi-agent settings. This class of games includes, for instance, any finite game and any multi-agent extension of traditional combinatorial optimization problems. After providing some background motivation and context of applications, we systematically review and classify the state-of-the-art algorithms to compute Nash equilibria. We propose an essential taxonomy of the algorithmic ingredients needed to compute equilibria, and we describe the theoretical and practical challenges associated with equilibria computation. Finally, we quantitatively and qualitatively compare a sequential Stackelberg game with a simultaneous \emph{IPG} to highlight the different properties of their solutions.
}
\iftoggle{ARXIV}{%
    \title{\theTitle}
    \date{} 
    \author{
        Margarida Carvalho \orcidA{} Gabriele Dragotto \orcidB{} Andrea Lodi \orcidC{}\\ Sriram Sankaranarayanan \orcidD{}
    }
                
    \maketitle     
    \begin{abstract}
        \theAbstract
    \end{abstract}
    }{

\CHAPTERNO{\phantom{Chapter 1}}%
\DOI{}
\TITLE{\theTitle}
\AUBLOCK{
    \AUTHOR{Margarida Carvalho}
    \AFF{Département d'informatique et de Recherche Opérationnelle, Universit\'e de Montr\'eal, Canada\\
        \EMAIL{carvalho@iro.umontreal.ca}}
    \AUTHOR{Gabriele Dragotto}
    \AFF{Department of Operations Research and Financial Engineering, Princeton University, U.S.A.
        \EMAIL{gdragotto@princeton.edu}}
    \AUTHOR{Andrea Lodi}
    \AFF{Jacobs Technion-Cornell Institute, Cornell Tech, and Technion - IIT, U.S.A.
        \EMAIL{andrea.lodi@cornell.edu}}
    \AUTHOR{Sriram Sankaranarayanan}
    \AFF{Operations and Data Sciences, IIM Ahmedabad, India\\
        \EMAIL{srirams@iima.ac.in}}
}
\CHAPTERHEAD{Carvalho et al. | Integer Programming Games: A Gentle Computational Overview}
\ABSTRACT{\theAbstract}

}
\maketitle

\section{Introduction}

How should an energy company determine its optimal production schedule? Which airlines should form a strategic alliance? How should graduate students from medical schools be matched with residency training? In all these contexts, decision-making is rarely an individual task; on the contrary, it often involves the mutual interaction of several self-interested stakeholders and their individual preferences. Neglecting the intrinsic role of strategic interaction may lead to inaccurate mathematical models and unreliable prescriptive recommendations. On top of the strategic interaction, decision-makers often face a concurrence of scarce resources (\eg, time, money, raw materials), complex operational constraints, and ample alternatives for their decisions.

In this tutorial, we provide an overview of \glspl{IPG}~\cite{koppe_rational_2011}, a large class of games embedding the modeling capabilities of mixed-integer programming. Specifically, these games combine the modeling power of integer variables with non-cooperative decision-making and provide a realistic and sophisticated tool to model, for instance, requirements such as indivisible quantities, logical disjunctions, and fixed production costs. While \glspl{IPG} include any finite game (\ie, games with a finite number of strategies and players), they enable us to implicitly represent any game whose set of outcomes is rather challenging to enumerate, \eg, the set of solutions of a difficult combinatorial problem. Therefore, \glspl{IPG} enable decision-makers to derive prescriptive recommendations that, on the one hand, account for complex operational requirements and, on the other hand, for strategic interaction with other decision-makers.

Compared to the solution of a single optimization problem, the solutions of \glspl{IPG} should account for the added complexity resulting from the strategic interaction. In this sense, solutions should embed the concept of \emph{stability}; namely, the solution should ensure that no player can unilaterally and profitably \emph{deviate} from it. In this context, the \emph{deviations} are context-dependent decisions, as they must reflect the realistic threats each decision-maker may face given its opponents' choices. For example, given the decisions of each energy company in an electricity market, an energy-producing firm may have the power to unilaterally decrease its energy production if the energy price is unattractive.
The leading solution concept embedding the notion of stability is arguably the Nash equilibrium~\cite{nash_equilibrium_1950,nash_noncoop_1951}. In a Nash equilibrium, all the decision-makers make mutually-optimal (\ie, rational) decisions, and no single decision-marker is incentivized to deviate to a different set of decisions. The ubiquitous concept of Nash equilibrium is one of the most fundamental concepts in game theory, with applications ranging from economics to social sciences. 

\vspace{4pt}

\paragraph{Tutorial Overview.} In this tutorial, we provide an overview of the motivations, the computational aspects, and the challenges associated with \glspl{IPG}. Specifically, we primarily focus on the algorithmic aspect connected to the computation of the Nash equilibria of \glspl{IPG}. In the remaining part of this section, we provide the mathematical formulation of \glspl{IPG} and some motivations for their use. In \cref{sec:litrev}, we provide a brief algorithmic literature review. In \cref{sec:blocks}, we propose a taxonomy to categorize the main building blocks employed by the available algorithms to compute equilibria and describe the algorithm's properties. In \cref{sec:criticalnode}, we present an example of an \gls{IPG} that models the protection of critical infrastructure and compare the solutions of this (simultaneous) game to the ones of a sequential version. Finally, in \cref{sec:ideas}, we provide ideas for future research and summarize this tutorial's key aspects.

\subsection{What}
\label{sec:what}
We can colloquially define an \gls{IPG} as a finite tuple of mixed-integer programming problems.
For each decision-maker, \ie, \emph{player}, we mathematically represent its space of possible decisions through a set of inequalities and integrality requirements, and its \emph{payoff} as a function of its decisions and those of the other players.
Let the operator $(\cdot{})^{-i}$ represent the elements of $(\cdot{})$ excluding the $i$-th component;  we formally define an \gls{IPG} in \cref{def:IPG}.
\begin{definition}{(\gls{IPG})}
    \label{def:IPG}
    An \gls{IPG} is a \emph{simultaneous}, \emph{non-cooperative} and \emph{complete-information} game among $n$ players where each player $i = 1,\dots,n $ solves the optimization problem
    \begin{subequations}
        \label{eq:IPG}
        \begin{align}
            \underset{\xofi}{\maxi} \quad & \{ f^i(\xofi;\xminusi) \text{ s.t. } x^i \in \X{i} \},    \text{ where }
                                           \X{i}:=\{x^i \in \mathbb{R}^{k_i}\times \mathbb{Z}^{n_i}: g^i(x^i)  \ge 0\}.
        \end{align}
    \end{subequations}
    Without loss of generality, we assume that: 
    \begin{enumerate*}
    \item $g^i: \mathbb{R}^{k_i+n_i} \rightarrow \mathbb{R}^{m_i}$, 
    \item each player controls $k_i$ continuous variables and $n_i$ integer variables, and
    \item each player has $m_i$ constraints defined by $g^i$.
    \end{enumerate*}
\end{definition}

We define any point $x^i$ feasible for the optimization problem of player $i$ as a \emph{pure strategy} for $i$. Accordingly, we define the feasible set $\X{i}$ as the \emph{strategy set} of player $i$, \ie, the set of decisions available for $i$. In many practical cases, and even according to the original definition of Köppe \etal~\cite{koppe_rational_2011}, the set $\X{i}$ is a mixed-integer set of the form $\{x^i \in \mathbb{R}^{k_i}\times \mathbb{Z}^{n_i} : A^i x^i \ge b^i\}$, where the entries of the matrix $A^i$ and the vector $b^i$ are rational numbers. We call the real-valued function $f^i(x^i;x^{-i})$ the \emph{payoff} function of player $i$, \ie, a function measuring the benefit (\eg, profit, satisfaction) that $i$ achieves by playing $x^i$ and observing $x^{-i}$ from its opponents. In \glspl{IPG}, the strategic interaction takes place in $f^i$, as the latter is a function of $x^i$ parametrized in $x^{-i}$, \ie, a function of the strategy of $i$ parametrized in the strategies of $i$'s opponents.
When we say the game is \emph{simultaneous}, we mean that each player commits to a strategy simultaneously with the other players, \ie, there is no order of play, and players cannot observe their opponents' strategies before committing to their strategies. When we say the game is \emph{non-cooperative}, we mean that players may have conflicting interests, for instance, conflicting payoff functions. When we say there is \emph{complete information}, we mean that every player is aware of its payoff function and strategy space and those of its opponents. Finally, when we say that each player is \emph{rational}, we mean that each player will maximize its payoff $f^i(x^i;x^{-i})$ and choose its strategy $x^i$ accordingly; for instance, given $\tilde{x}=(\tilde{x}^1,\dots,\tilde{x}^n)$, player $i$ would never play any $\hat{x}^i$ such that $f^i(\hat{x}^i;\tilde{x}^{-i}) <  f^i(\tilde{x}^i;\tilde{x}^{-i})$.
We consider the Nash equilibrium of \cref{def:NE} as the solution concept for the \gls{IPG} of \cref{def:IPG}.
\begin{definition}[Pure Nash equilibrium]
    A \emph{pure Nash equilibrium} for an \gls{IPG} is a vector of pure strategies $\bar{x}=(\bar{x}^1,\dots,\bar{x}^n)$ such that, for each player $i=1,\dots,n$,
    \begin{align}
        f^i(\bar{x}^i;\bar{x}^{-i}) \ge f^i(\tilde{x}^i;\bar{x}^{-i}) \quad \forall \tilde{x}^i \in \X{i}.
        \label{eq:NE}
    \end{align}
    \label{def:NE}
\end{definition}
\paragraph{Mixed and Approximate Equilibria.} 
Pure equilibria are also called \emph{deterministic} since each player commits to a strategy with unit probability. However, \cref{def:NE} also extends to \emph{mixed} strategies and equilibria, \ie, randomized strategies and equilibria. Specifically, whenever players randomize over their strategies, instead of playing a pure strategy deterministically, we call the strategy \emph{mixed}. A mixed strategy is a probability distribution over the set of pure strategies; we will let $\Delta^i$ be the \emph{mixed strategy set} of player $i$. A mixed Nash equilibrium is then a vector of strategies where \cref{eq:NE} holds over the expectations of the players' payoff, as in \cref{def:MNE}.
\begin{definition}[Mixed Nash equilibrium]
    A \emph{mixed Nash equilibrium} for an \gls{IPG} is a vector of mixed strategies $\bar{x}=(\bar{x}^1,\dots,\bar{x}^n)$ such that, for each player $i=1,\dots,n$, $\bar{x}^i \in \Delta^i$ and
    \begin{align}
        \mathbb{E}_{X \sim \bar{x}} \big(f^i(X^i;X^{-i}) \big) \ge \mathbb{E}_{X \sim \bar{x}}\big( f^i(\tilde{x}^i;X^{-i}) \big) \quad \forall \tilde{x}^i \in \X{i}.
        \label{eq:MNE}
    \end{align}
    \label{def:MNE}
\end{definition}
We remark that, by the definition of Nash equilibrium, in \cref{def:MNE} we restrict the search for possible \emph{deviations} $\tilde{x}^i$ to the pure strategies. 
Finally, whenever \cref{eq:NE} (\cref{eq:MNE}) is approximately satisfied, \eg, with a relative or absolute violation of a small positive quantity, the equilibrium is an \emph{approximate} Nash equilibrium. We conclude by showcasing a simple \gls{IPG} and its equilibria in \cref{ex:knapsack}.

\begin{example}[Knapsack Game]
\label{ex:knapsack}
Consider a 2-player \gls{IPG} where each player solves a binary knapsack problem as in \cref{eq:knapsack}. Specifically, the players solve the parametrized integer programs
\begin{subequations}
\begin{align}
\underset{x^1 \in \{0,1\}^2}{\text{maximize}} \big\{ & x^1_1 + 2x^1_2 - 2x^1_1x^2_1 -3x^1_2x^2_2 & \!\!\!\!\!\!\!\!\!\!\!\!\!\!\!\!\!\!\!\!\!\!\!\!\!\!\!\!\!\!\!\!\!\!\!\!\!\!\text{ s.t. }\;\;& 3x^1_1+4x^1_2 \le 5 \big\}, \\
\underset{x^2 \in \{0,1\}^2}{\text{maximize}}  \big\{ &  3x^2_1 + 5x^2_2 - 5x^2_1x^1_1 -4x^2_2x^1_2 & \!\!\!\!\!\!\!\!\!\!\!\!\!\!\!\!\!\!\!\!\!\!\!\!\!\!\!\!\!\!\!\!\!\!\!\!\!\!\text{ s.t. }\;\;&
2x^1_1+5x^1_2 \le 5 \big\}.
\end{align}
\label{eq:knapsack}
\end{subequations}
Player $i$ can select a given \emph{item} $j=1,2$ by setting $x^i_j=1$ and receive a profit that solely depends on its choice; for instance, by selecting item $2$, player $2$ gets a linear profit of $5$. However, when both players select a given item, the interaction alters their payoffs; for instance, when they both select item $2$, the payoff of player $2$ gets penalized by $-4$. Because of their packing constraints, the only feasible strategies for the players are $(x^i_1,x^i_2) \in \{ (0,0), (0,1), (1,0)\}$. We can easily verify that this \gls{IPG} has 3 Nash equilibria: the pure Nash equilibria
$(x^1_1,x^1_2,x^2_1,x^2_2) = (0, 1, 1, 0)$ and $(1, 0, 0, 1)$, and the mixed Nash equilibrium $(2/9, 7/9, 2/5, 3/5)$. In the latter, player $1$ plays the strategy $(x^1_1,x^1_2)=(1,0)$ with probability $2/9$ and the strategy $(x^1_1,x^1_2)=(0,1)$ with probability $7/9$; symmetrically, player $2$ plays the strategy $(x^2_1,x^2_2)=(1,0)$ with probability $2/5$ and the strategy $(x^2_1,x^2_2)=(0,1)$ with probability $3/5$.
\end{example}

\subsection{Why}
\label{sub:why}
There are three substantial reasons behind the interest in \glspl{IPG}: their modeling capabilities, the informative content of their Nash equilibria, and the growing empirical evidence of their practical use cases.

\vspace{4pt}

\paragraph{Modeling Capabilities.}
As previously mentioned, \glspl{IPG} combine the modeling power of mixed-integer programming with the strategic interactions stemming from a game setting. Although the integrality requirements on the variables make the optimization problem more challenging, the theoretical and practical advances in integer programming propelled dramatic performance improvements in mixed-integer programming. Indeed, in the absence of strategic interactions, we can nowadays reliably and efficiently solve large-scale mixed-integer optimization problems with open-source and commercial solvers. \glspl{IPG} inherit the modeling capabilities, as well as the difficulty, of mixed-integer programming and complement the latter by enabling the modeling of strategic interactions. Formulating an \gls{IPG} requires an implicit description of the player's strategies, \eg, in terms of a set of inequalities and integrality requirements, and of the payoff functions. While any simultaneous game is a tuple of optimization problems, the core difference between classic representations (\eg, normal-form or bimatrix games) and \glspl{IPG} lies precisely in the implicit description of the players' strategy sets. For instance, if the set of strategies for a given player coincides with the set of solutions of a difficult combinatorial problem, the representation of the game in a classic form would require the enumeration of the exponentially many feasible points of the combinatorial problem.
On the contrary, \glspl{IPG} can compactly represent games without explicitly enumerating the players' strategies. Finally, \glspl{IPG} are inherently nonconvex games due to the integrality requirements on some of the players' variables. While plenty of literature focuses on convex games, little is known about games with non-convexities, both in theory and practice~\cite{daskalakis_nonconvex_2022}. \glspl{IPG} are no exception, and, despite their modeling power, the methodological and practical research is relatively recent.

\vspace{4pt}

\paragraph{Nash Equilibria.} From a qualitative perspective, Nash equilibria are stable solutions in the sense that no player can profitably deviate from them. As such, game theorists often assume it is reasonable to expect that players will play Nash equilibria under the assumptions of complete information and rationality. However, from an external observer's perspective, the concept of Nash equilibrium is even more pivotal, as it enables the development of regulatory interventions. Indeed, Nash equilibria are solutions trading off the individual selfishness of each player and overarching societal goals, for instance, fairness or social welfare. Alvin Roth~\cite{roth_economist_2002}, with his argument “economists as engineers”, strongly motivated the need for economists and regulators to not only analyze strategic interactions (\eg, markets) but rather to \emph{design} them as engineers would do. In this context, design means determining the best market rules (\eg, incentives, taxes, penalties) so that the market exhibits some desirable properties. Nash equilibria provide, in this sense, the endpoint of these analyses, while \glspl{IPG} can provide the complex modeling tools the practitioners need. Therefore, computing equilibria is not only a methodological exercise but a fundamental step for developing sophisticated economic and mathematical models and regulatory policies. As \glspl{IPG} are games without specific structure (besides the few assumptions of \cref{def:IPG}), the algorithms for computing their equilibria make broad and unrestrictive assumptions on the game's structure instead of the limiting and punctual assumptions economists may have to make on their models to ensure the computability or uniqueness of equilibria.

\vspace{4pt}

\paragraph{Practical Use Cases.} Although \glspl{IPG} are a relatively new class of games, practitioners and researchers have already employed them in some applied settings, for instance, revenue management, healthcare, cybersecurity, and sustainability. We provide a non-exhaustive review of some of these use cases to illustrate the flexibility of \glspl{IPG} as modeling tools. In production economics, Lamas and  Chevalier~\cite{LAMAS2018864} model a game where players make pricing and inventory-management decisions; in this context, each player solves a mixed-integer optimization problem adapted from the well-known mixed-integer formulations of the lot-sizing problem. Similarly, Carvalho \etal~\cite{carvalho2018} let players decide on the quantity of goods instead of prices.  Crönert and Minner~\cite{cronert_equilibrium_2020} model a facility location and design game where the players make strategic decisions on where to geographically locate facilities and what assortment to use to maximize profits. Sagratella \etal~\cite{sagratella_noncooperative_2020} employ \glspl{IPG} to provide a game-theoretic extension of the fixed-charge transportation problem.
In cybersecurity, Dragotto \etal~\cite{Dragotto_critical_2023} develop an \gls{IPG} model to assess the cyber-security risk of cloud networks and inform security experts on the optimal security strategies to adapt in the eventuality of a cyber-attack. In healthcare,~\cite{Blom2022,carvalho_nash_2017} model cross-border kidney exchange programs via \glspl{IPG}. In this latter applicative context, the algorithms for solving \glspl{IPG} enable regulators to evaluate different exchange-mechanism designs for those kidney exchange programs; namely, they help determine whether the program's mechanism (rules) leads to equilibria maximizing the patients' benefits (\eg, equilibria maximizing the number of kidney transplants).
Finally, in energy, Crönert and Minner~\cite{cronert_2021} model the competition between providers of hydrogen fuel stations as an \gls{IPG}.

\section{Computing Equilibria: Challenges and Algorithms}
\label{sec:litrev}

This section provides a brief literature review of the available methodologies to compute Nash equilibria in \glspl{IPG}. In \cref{sec:blocks}, we will provide a more detailed description of the main algorithms and an analysis of their main components. Although a wealth of theories and algorithms exist to compute Nash equilibria in well-structured (\eg, normal-form and extensive-form) games, the literature on \glspl{IPG} is significantly more recent and limited. 

\vspace{4pt}

\paragraph{Two Limitations.} Most approaches for computing Nash equilibria in games focus on normal-form or convex games.
On the one hand, several foundational works (\eg, see~\cite{lemke_equilibrium_1964,rosenmuller_generalization_1971,sandholm_mixed-integer_2005}) consider finite games represented in normal form, \ie, through matrices describing the payoffs under any combination of the players' strategies. As previously mentioned, this representation implicitly requires the enumeration of the players' strategies, which is impractical in \glspl{IPG}; indeed, the respective normal-form representation may have an exponential size or even be impossible to characterize due to unbounded or uncountable strategy sets.
On the other hand, whenever the game is continuous, \ie, for each player, the set of strategies is a non-empty compact metric space, and the payoff is a continuous function, the task of computing equilibria is equivalent, under some technical assumption, to solving a variational inequality or complementarity problem~\cite{facchinei_finite-dimensional_2003}. Although this class of algorithms is efficient in practice, it does not handle the non-convexities induced by the integer variables.

\vspace{1em}
While the above limitations prevent the direct application of normal-form and variational or complementarity-based techniques to the computation of equilibria in \glspl{IPG}, they also inspired several authors to extend these techniques with new algorithmic components. This tutorial will analyze and characterize the six algorithms, which we briefly introduce in the following paragraphs; we will later analyze them in detail in \cref{sec:blocks}.

\vspace{4pt}

\paragraph{Mixed Equilibria.} We analyze three algorithms to compute mixed equilibria in \glspl{IPG}. First, we consider the \gls{SGM} from Carvalho \etal~\cite{carvalho_computing_2017}, an algorithm to compute a mixed Nash (and correlated) equilibrium in \glspl{IPG} with separable and Lipschitz continuous payoff functions; separable payoffs take a sum-of-products form~\cite{stein_separable_2008}. 
Second, we analyze the \gls{eSGM} from Crönert and Minner~\cite{cronert_equilibrium_2020}, an algorithm extending \gls{SGM} that can enumerate all the mixed Nash equilibria when players only control integer variables (\ie, $k_i=0$). Finally, we analyze \gls{CNP} from Carvalho \etal~\cite{Dragotto_2021_CNP}, an algorithm to compute an equilibrium for nonconvex games with separable payoffs, including \glspl{IPG}, by outer approximation of the set of strategies for each player; the algorithm exploits an equivalent convex reformulation of the game where each player $i$ strategy set is $\cl \conv(\X{i})$, \ie, the closure of the convex hull of its original strategy set, as opposed to $\X{i}$.

\vspace{4pt}

\paragraph{Pure Equilibria.} Several authors focused on the computation of pure Nash equilibria~\cite{sagratella_computing_2016,sagratella_2019,ZERORegets,schwarze_branch-and-prune_2022,harks_generalized_2022,guo_copositive_2021}. In particular, the works in~\cite{sagratella_2019,ZERORegets,harks_generalized_2022} allow the players' constraints to depend on the opponents' variables, hence computing the so-called \emph{generalized} pure Nash equilibria. For the purpose of this tutorial, we will focus on three algorithms to compute pure Nash equilibria and omit the generalized equilibria; we will provide further pointers to generalized equilibria in \cref{sec:ideas}.
First, we consider the \gls{BM} from Sagratella~\cite{sagratella_computing_2016}, an algorithm to compute and enumerate pure equilibria whenever the players have continuously-differentiable and concave payoffs in their variables and feasible sets described via convex inequalities and integrality requirements. 
Second, we analyze the \gls{BNP} method from Schwarze and Stein~\cite{schwarze_branch-and-prune_2022}, an algorithm that extends \gls{BM} by requiring less restrictive assumptions on the players' payoff functions.
Finally, we analyze \gls{ZEROR} from Dragotto and Scatamacchia~\cite{ZERORegets}, a cutting plane algorithm to optimize over the set of pure Nash equilibria of \glspl{IPG} with mixed-integer linear strategy sets.

\vspace{0.8em}

As the literature on \glspl{IPG} is relatively recent, we observe no dominant algorithm among the presented ones. Furthermore, we remark that the algorithms often use different algorithmic components and support different types of \glspl{IPG} (\eg, different payoff functions). The following section provides an abstract taxonomy of the main algorithmic components of the six algorithms we introduced to provide a unified framework for their analysis.

\section{Solving IPGs}
\label{sec:blocks}

From a practical perspective, the plausibility of \glspl{IPG} as modeling tools also stems from the availability of efficient tools to compute their equilibria. However, even besides the realm of \glspl{IPG}, the computation of Nash equilibria requires researchers to address some fundamental challenges. This section provides an abstraction of the main \emph{building blocks} (\ie, components) employed in the available algorithms and analyzes the practical challenges associated with the computation of Nash equilibria.

\subsection{Algorithmic Building Blocks}
\label{sub:blocks}
In order to embrace the diversity of the available algorithms to compute Nash equilibria in \glspl{IPG}, we propose an abstraction of their main components. We will later contextualize these components for each algorithm.

\begin{figure}[!ht]
    \centering
    \includegraphics[width=0.95\textwidth]{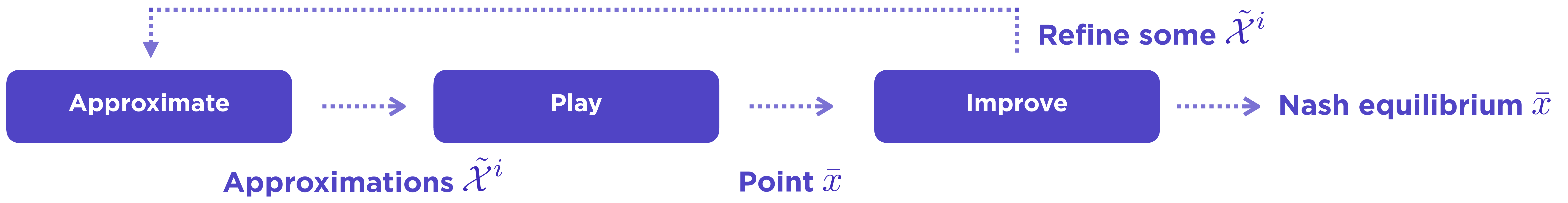}
    \caption{The three-phase process to compute Nash equilibria in \glspl{IPG}.}
    \label{fig:phases}
\end{figure}

We identify three main \emph{building blocks} shared by the majority of the algorithms we review: an \emph{approximation} phase, a \emph{play} phase, and an \emph{improvement} phase. In the interest of clarity and without loss of generality, we will assume that at least one Nash equilibrium exists, and we will present our analysis by referring to \emph{the algorithm} as an abstract routine. We schematically represent the three-phase process in \cref{fig:phases} and in \cref{Alg:Abstract}.

\begin{algorithm}[!ht]
    \SetKwBlock{Repeat}{repeat}{}
    \DontPrintSemicolon
    \caption{An abstract algorithm to compute Nash equilibria for \glspl{IPG} \label{Alg:Abstract}}
    \KwData{An \glspl{IPG} instance.}
    \KwResult{A Nash equilibrium $\bar{x}$.}
    \textbf{(Approximate)} Determine the initial approximation $\mathcal{\tilde{X}}^i$ for each player $i$\;
    \Repeat
    {
        \textbf{(Play)} Solve an optimization problem to determine a tentative equilibrium $\bar{x}$\;
        \textbf{(Improve)} $A\gets$ the answer from the improvement oracle (\cref{Alg:AbstractOracle})\;
        \lIf{$A$ is \texttt{yes}}{
            \KwRet{the Nash equilibrium $\bar x$}
        }
        \Else{
            \textbf{(Approximate)} Refine $\mathcal{\tilde{X}}^i$ for some player $i$ such that $\bar x$ is infeasible.
        }
    }
\end{algorithm}

\paragraph{Approximate.} In the approximation phase, the algorithm selects an approximated strategy space $\tilde{\mathcal{X}}^i$ for each player $i$. 
Typically, these approximated strategy spaces exhibit desirable properties (\eg, convexity), and identifying equilibria for these approximations is relatively easier than for the original strategy spaces; furthermore, the approximated strategy space either enlarges or restricts the space of available strategies for each player.
Whenever the approximation $\tilde{\mathcal{X}}^i$ contains the original strategy space $\X{i}$, we call it an \emph{outer approximation}, \ie, a set $\tilde{\mathcal{X}}^i \supseteq \X{i}$ that possibly introduces infeasible (with respect to $\X{i}$) strategies for player $i$; otherwise, if the approximation $\tilde{\mathcal{X}}^i$ possibly removes some strategies belonging to $\tilde{\mathcal{X}}^i$, we call it an \emph{inner approximation}, \ie, $\tilde{\mathcal{X}}^i \subseteq \X{i}$.

\vspace{4pt}

\paragraph{Play.} In the \emph{play} phase, the algorithm usually computes a tentative Nash equilibrium $\bar{x}=(\bar{x}^1,\dots,\bar{x}^n$) in the \emph{approximated} game originating from the approximated strategy spaces selected in the previous phase, \ie, the game where each player $i$'s strategy space is $\tilde{\mathcal{X}}^i$. The play phase often involves solving an optimization problem whose particular structure heavily depends on the specific algorithm. The algorithms we review in this paper often solve complementarity, linear, or mixed-integer optimization problems. For instance, \gls{CNP} solves a complementarity problem, \gls{ZEROR} solves a mixed-integer program, and \gls{SGM} solves a series of (feasibility) problems, typically linear (\cref{ex:play}).  The optimization problem solved in this phase is generally formulated over a feasible region that depends on all the players' variables, \ie, it depends on $x^1,\dots,x^n$.

\begin{example}[Play Phase of \gls{SGM}]
    \label{ex:play}
    Consider an \gls{IPG} and let, for each player $i$, $\tilde{\X{i}}$ be a subset of the points in $\X{i}$. In the play phase, \gls{SGM} solves a set of feasibility programs which are linear if $n=2$, or if the interactions between opponents' variables are bilinear, \ie, the terms in $f^i$ accounting for the opponents' variables are in the form $(x^i)^\top Q^i x^{-i}$ where $Q^i$ is a matrix of appropriate dimensions. In the play phase, for each choice of support $S^i \subseteq \tilde{\X{i}}$ and player $i$, \gls{SGM} seeks to find a mixed equilibrium $\sigma=(\sigma^1,\dots,\sigma^n)$ with a vector of expected payoffs $v=(v^1,\ldots,v^n)$ by solving the following feasibility problem
        \begin{align*}
            \sigma^i_j \ge 0 \qquad               & \forall i=1,\ldots,n,\forall j \in S^i,                           \\
            \sigma^i_j = 0 \qquad                 & \forall i=1,\ldots,n,\forall j \in \tilde{\X{i}} \setminus S^i,   \\
            \sum_{j \in S^i} \sigma^i_j =1 \qquad & \forall i=1,\ldots,n,                                             \\
            f^i(x^i,\sigma^{-i}) = v^i \qquad     & \forall i=1,\ldots,n, \forall x^i \in S^i,                        \\
            f^i(x^i,\sigma^{-i}) \leq v^i\qquad   & \forall i=1,\ldots,n,\forall x^i \in \tilde{\X{i}} \setminus S^i.
        \end{align*}
        The first three constraints enforce that $\sigma^i$ is a probability distribution over $S^i$, while the last two constraints guarantee that the strategies in the support $S^i$ are best responses.
\end{example}

\paragraph{Improve.} Given the tentative equilibrium $\bar{x}$, the improvement phase determines whether $\bar{x}$ is also a Nash equilibrium for the original game and otherwise refines the approximation, \eg, by excluding $\bar{x}$ from the approximations of the first phase, or by adding some new strategies to the approximations. 
We define the routine performing this improvement phase as the \emph{improvement oracle}, which we sketch in \cref{Alg:AbstractOracle}.
Although different algorithms implement different improvement oracles, most of the improvement oracles involve two steps.
First, the oracle performs a \emph{membership} step, where the routine checks if, for any player $i$, $\bar{x}^i$ is a feasible strategy for the original game; for instance, the oracle could check if $\bar{x}^i \in \X{i}$ when the algorithm computes pure equilibria; more generally, when the algorithm computes mixed equilibria, the oracle could check if the strategies in the support of $\bar{x}^i$ are feasible pure strategies for the original game.
Naturally, the membership step is trivial when the approximation of the first phase is an inner approximation. 
Second, the improvement oracle performs a \emph{stability} step, where the routine checks, for any player $i$, whether the tentative strategy $\bar{x}^i$ is a \emph{best response} for player $i$ to $\bar{x}^{-i}$; namely, it checks if $\bar{x}^i$ is the best strategy that $i$ can play in response to $\bar{x}^{-i}$, or whether there exists a deviation from $\bar{x}^i$ guaranteeing player $i$ a better payoff. In this step, the routine determines the optimal solution $\tilde x^i$ to the player's optimization problem by fixing $x^{-i}$ to the tentative strategy $\bar x^{-i}$, \ie, it solves the optimization problem $\max_{x^i} \{ f^i(x^i;\bar x^{-i}):x^i \in \X{i} \}$. If the payoff guaranteed by $\tilde x^i$ is greater than the one guaranteed by $\bar x^i$ under $\bar{x}^{-i}$, then $\bar x^i$ can not be a best response, \ie, if $f^i(\tilde x^i;\bar{x}^{-i}) >
    f^i(\bar x^i;\bar{x}^{-i})$.
Suppose the improvement oracle determines, for every player $i$, that $\bar{x}^i$ is both feasible and a best response. In that case, it terminates by providing a \texttt{yes} answer, \ie, certifying that $\bar x$ is a Nash equilibrium for the \gls{IPG}. Otherwise, the routine returns a \texttt{no} and some additional information $\mathbf{I}$, \eg, a best response $\tilde x^i$ for some player $i$.
Starting from $\mathbf{I}$, the algorithm refines the successive approximations and iterates again the play phase. In an outer approximation algorithm, $\mathbf{I}$ can contain separating hyperplanes that refine the players' approximations, \eg, see \gls{ZEROR}, \gls{BM} and \gls{CNP}. We showcase how to employ $\mathbf{I}$ in \cref{ex:improve_SGM,ex:improve_CNP}. We remark that, whenever the algorithm computes mixed Nash equilibria, we should be able to interpret the solution of the play phase as a mixed Nash equilibrium of the original game, \ie, we should have access to a routine that transforms the mixed strategy of each player to a probability distribution over the strategies of the approximated game. Equivalently, we should be able to determine the strategies played with positive probabilities by each player.

\begin{example}[Improve Phase of \gls{SGM}]
    \label{ex:improve_SGM}
    Consider an \gls{IPG} and let, for each player $i$, $\tilde{\X{i}}$ be a subset of the points in $\X{i}$ at a given iteration of the algorithm. The improvement oracle, since $\tilde{\X{i}} \subseteq \X{i}$, does not need to perform the membership step, as it is trivially satisfied.
    Assume that, in the stability step, there exists a player $i$ that can deviate to $\tilde{x}^i$ and improve its payoff. Then, the information $\mathbf{I}$ contains $\tilde{x}^i$, and the algorithm can refine $\tilde{\X{i}}$ by including $\tilde{x}^i$, \ie, $\tilde{\X{i}}=\tilde{\X{i}}\cup\{\tilde{x}^i\}$.
\end{example}

\begin{algorithm}[!ht]
    \SetKwBlock{Repeat}{repeat}{}
    \DontPrintSemicolon
    \caption{An abstract improvement oracle for \glspl{IPG} \label{Alg:AbstractOracle}}
    \KwData{An \glspl{IPG} instance and a tentative equilibrium $\bar x$.}
    \KwResult{Either \begin{enumerate*} \item \texttt{yes}, or \item \texttt{no} and some additional information $\mathbf{I}$. \end{enumerate*}}
    \ForEach{player $i = 1,\dots,n$}{
        \If{\textbf{(Membership)} $\bar x^i \in \X{i}$ (or, any strategy in the support of $\bar x^i$ is in $\X{i}$)}{
            $\tilde x^i \gets \arg\max_{x^i}\{ f^i(x^i;\bar x^{-i}) : x^i \in \X{i} \}$\;
            \If{\textbf{(Stability)} $f^i(\bar x^i; \bar{x}^{-i}) \ge f^i(\tilde x^i; \bar{x}^{-i})$ }{
                $\bar x^i$ is a feasible best response for player $i$
            }
            \Else{
                \KwRet{\texttt{no} and $\mathbf{I}$ (\eg, $\tilde x^i$)}
            }
        }
        \Else{
            \KwRet{\texttt{no} and $\mathbf{I}$}
        }
    }
    \KwRet{\texttt{yes}}
\end{algorithm}

\begin{example}[Improve Phase of \gls{CNP}]
    \label{ex:improve_CNP}
Consider an \gls{IPG} and let, for each player $i$, $\tilde{\X{i}}$ be a superset of $\X{i}$ at a given iteration of \gls{CNP}, \ie, $\tilde{\X{i}} \supseteq \X{i}$. In contrast to \cref{ex:improve_SGM}, the improvement oracle needs to perform the membership step to determine if $\bar{x}^i \in \tilde{\X{i}}$, or, for mixed Nash equilibria, if any strategy in the support of $\bar{x}^i$ belongs to $\tilde{\X{i}}$.
    In \gls{CNP}, Carvalho \etal~\cite{Dragotto_2021_CNP} prove that this membership step is equivalent to determining if $\bar{x}^i \in \cl \conv(\X{i})$ for any \gls{IPG} where the payoff functions are separable. If $\cl \conv (\X{i})$ is a polyhedron (\ie, the so-called \emph{perfect formulation} of a mixed-integer programming set), this is equivalent to solving the following \emph{separation task}: given $\bar{x}^i$ and $\X{i}$, determine if $\bar{x}^i \in \cl \conv(\X{i})$ or provide a separating hyperplane $\pi^\top x^i \ge \pi_0$ such that $\pi^\top x^i \ge \pi_0$ for any $x^i \in \X{i}$, and $\pi^\top \bar{x}^i < \pi_0$. If $\bar{x}^i \notin \X{i}$, the improvement oracle will return \texttt{no} and $\mathbf{I}$ will contain a separating hyperplane $\pi^\top x^i \ge \pi_0$; in practice, the latter can refine the approximation $\tilde{\X{i}}$.
    Solving the separation task and, generally, refining the approximations, involves traditional integer-programming methodologies such as branching and cutting plane generation \cite{conforti_integer_2014}. 
\end{example}

\subsection{Theoretical and Practical Challenges}
Despite the efforts to develop efficient methods to compute Nash equilibria, the task still involves several theoretical and computational challenges. 

\vspace{4pt}

\paragraph{Choosing the Approximation.} In contrast to the solution of approximated (\ie, relaxed, inner, and outer approximated) optimization problems, the Nash equilibrium of a game's approximation does not imply the existence of equilibria in the original game. We illustrate this in \cref{ex:approximations}.

\begin{example}[Hierarchy of Approximations in \glspl{IPG}]
    \label{ex:approximations}
    Consider an \gls{IPG} with two players such that they solve
    \begin{align*}
        \underset{x^1}{\maxi} \{ -x^1 x^2 : x^1 \ge 1, x^1 \in \mathbb{Z} \}, \qquad \qquad
        \underset{x^2}{\maxi} \{ x^2x^1 : x^2 \in \{-1,1\} \}.
    \end{align*}
    This game has one Nash equilibrium, namely $x^1=1$ and $x^2=1$. However, different inner or outer approximations of the players' strategy sets may result in approximated games exhibiting different equilibria compared to the original game. \\
    \emph{Approximation 1.} Let $\tilde{\X{1}}$ be $\X{1} \cap \{ x^1 \ge 2\}$, and $\tilde{\X{2}}$ be $\X{2}$; this approximated game has an equilibrium $x^1=2$ and $x^2=1$, which is not an equilibrium for the original game.\\
    \emph{Approximation 2.} Let $\tilde{\X{2}}$ be $\{ x^2 \ge -1, x^2 \in \mathbb{Z}\}$, and $\tilde{\X{1}}$ be $\X{1}$; this approximated game does not admit any equilibrium.
\end{example}
Indeed, the relationship between the Nash equilibria of an \gls{IPG} and those of its approximations intrinsically influences the algorithmic design of the three phases. For instance, selecting an unbounded approximation $\tilde{\X{i}}$ may result in the failure of the play phase to provide a tentative equilibrium $\bar x$, as an equilibrium might not exist.

\vspace{4pt}

\paragraph{Existence.}
The task of determining if an equilibrium exists may turn out to be difficult. While Nash proved that an equilibrium always exists in finite games and Glicksberg~\cite{glicksberg_further_1952} proved that an equilibrium always exists in continuous games, an equilibrium might not exist in an \gls{IPG}.  Based on these results, Carvalho \etal~\cite{vaz_existence_2018} stated that a sufficient condition for the existence of equilibria in \glspl{IPG} is the boundedness of the $\mathcal{X}^i$. If some strategy space is unbounded,~\cite{vaz_existence_2018,carvalho_computing_2017} established that the task of determining if an \gls{IPG} admits a Nash equilibrium (even if restricted to pure equilibria) is a \SigmaTwoP problem. Problems in the \SigmaTwoP class are in the so-called “second level of the polynomial hierarchy”, and unless the hierarchy collapses, it is impossible to represent them as integer programs of polynomial size.

\vspace{4pt}

\paragraph{Computing Tentative Equilibria.} In practice, the existence of equilibria intrinsically affects the design of algorithms since the latter may need to either deal with finite \glspl{IPG} or detect the non-existence (\eg, infeasibility) of an equilibrium. First, the issue of existence affects the play phase, as it requires the algorithm to compute the tentative equilibrium or determine its non-existence efficiently. Second, it also affects the approximation phase, as the algorithm should be able to iteratively refine the approximations depending on the output of the play phase. Moreover, even if the \gls{IPG} is finite, or it has a non-empty bounded set of strategies, and thus a (mixed) Nash equilibrium always exists, computing it is computationally demanding; indeed, the play phase often involves solving an \NPhard optimization problem, \eg, a complementarity problem or a mixed-integer optimization problem. From a theoretical perspective, this is not surprising given that determining the existence of Nash equilibria in \glspl{IPG} is \SigmaTwoP and, even in normal-form games (a particular case of \glspl{IPG}), the task of computing an equilibrium is difficult (\ie, it is PPAD-complete \cite{papadimitriou1994complexity}). From a practical perspective, however, this implies that the algorithms for \glspl{IPG} may require double-exponentially many iterations (with respect to the binary encoding of the instance) to converge to an equilibrium.

\vspace{4pt}

\paragraph{Selection.} The existence of multiple Nash equilibria exposes the issue of selection, \ie, the task of selecting a specific Nash equilibrium or, from an optimization perspective, the task of optimizing an objective function over the set of Nash equilibria. This is a different definition from the one of Harshany~\cite{harsanyi_new_1995} or  Harshany and Selten~\cite{harsanyi_new_1995}, as it does not directly deal with the \emph{probability} of observing an equilibrium (or with the notion of stability of an equilibrium) but rather of selecting an equilibrium based on a given function of the players' variables. For instance, the objective function can model concepts such as social welfare, fairness, or societal return. As the research on computing equilibria in \glspl{IPG} is relatively recent, only some algorithms are currently capable of \emph{selecting} equilibria, as opposed to computing \emph{an} equilibrium. We also argue that there is a fundamental difference between the selection and the enumeration of equilibria, as there is a difference between finding the optimal solution to an optimization problem and enumerating its solutions.

\subsection{Classifying the Algorithms}
Following the abstract three-phase process introduced in \cref{sub:blocks}, we classify the algorithms of \cref{sec:litrev} in two families based on the approximation phase: \emph{inner} and \emph{outer} approximation algorithms. As the name suggests, the main difference stems from the choice of approximation. In \cref{tab:classification}, we summarize the available algorithms for \glspl{IPG} and their characteristics.

\begin{table}[!ht]
    \centering
    \resizebox{\textwidth}{!}{
        \begin{tabular}{l@{\hspace{7em}}c@{\hspace{3em}}c@{\hspace{3em}}c@{\hspace{2em}}c}
            \hline
            \textbf{Method}                                 & \textbf{Class} & \textbf{Type of Equilibria} & \textbf{Equilibria Selection} & \textbf{Equilibria Enumeration} \\
            \hline
            \\
            \gls{SGM}~\cite{carvalho_computing_2017}        & Inner          & Mixed                       & \xmark                        & \xmark                          \\
            \gls{eSGM}~\cite{cronert_equilibrium_2020}      & Inner          & Mixed                       & \cmark                        & \cmark                          \\
            \gls{CNP}~\cite{Dragotto_2021_CNP}              & Outer          &
            Mixed                                           & \xmark         & \xmark                                                                                        \\
            \gls{BM}~\cite{sagratella_computing_2016}       & Outer          & Pure                        & \xmark                        & \cmark                          \\
            %
            %
            %
            %
            \gls{ZEROR}~\cite{ZERORegets}                   & Outer          & Pure$^\dagger$              & \cmark                        & \cmark                          \\
            \gls{BNP}~\cite{schwarze_branch-and-prune_2022} & Outer          & Pure                        & \xmark                        & \cmark
            \\
            \hline
        \end{tabular}
    }
    \vspace{2em}
    \resizebox{\textwidth}{!}{
        \begin{tabular}{l@{\hspace{1em}}c@{\hspace{1em}}c@{\hspace{1em}}c}
            \hline
                                                            & \textbf{Strategy Sets} $\bm{\X{i}}$ & \textbf{Payoffs} $\bm{f^i}$ & \textbf{Play Phase}     \\
            \hline
            \\
            \gls{SGM}~\cite{carvalho_computing_2017}        & Bounded mixed-integer linear                & Separable and Lipschitz continuous in $x^i$                   & Feasibility problems          \\
            \gls{eSGM}~\cite{cronert_equilibrium_2020}      & Bounded integer and inequality constraints                      & No restriction                   & Mixed-integer problem   \\
            \gls{CNP}~\cite{Dragotto_2021_CNP}              & Polyhedral convex hull $\cl\conv(\X{i})$        & Linear in $x^i$             & Complementarity problem \\
            \gls{BM}~\cite{sagratella_computing_2016}       & Bounded convex-integer       & Concave in $x^i$             & Complementarity problem   \\
            \gls{ZEROR}~\cite{ZERORegets}                   & Bounded mixed-integer linear$^\dagger$               & MIP linearizable            & Mixed-integer problem   \\
            \gls{BNP}~\cite{schwarze_branch-and-prune_2022} & Bounded convex-integer       & Strictly
concave in $x^i_j$       & Mixed-integer problem 
            \\
            \hline
        \end{tabular}
    }
    \caption{Summary of the available algorithms to compute Nash equilibria for \glspl{IPG}.
        \\ $^\dagger$: $\X{i}$ can include linearizable constraints on $x^{-i}$.}
    \label{tab:classification}
\end{table}

\begin{enumerate}
    \item \gls{SGM}~\cite{carvalho_computing_2017} is an inner approximation algorithm to compute a Nash equilibrium in \glspl{IPG} with separable payoff function and bounded mixed-integer linear strategy spaces. In its approximation phase, the algorithm \emph{samples} a restricted subset of pure strategies for each player and builds an inner approximation of the original game. In the play phase, the algorithm represents the inner-approximated game as a normal-form game and finds an equilibrium with the method of Porter \etal~\cite{porter_simple_2008} (recall \cref{ex:play}). In the improvement phase, the algorithm only performs the stability step, as the tentative equilibrium trivially passes the membership step. Whenever the stability steps determine a deviation, the algorithm updates its inner approximations and repeats the play phase. The algorithm does not support equilibria enumeration or equilibrium selection. Finally, as the authors show, \gls{SGM} can also compute correlated equilibria.
    \item \gls{eSGM}~\cite{cronert_equilibrium_2020} refines \gls{SGM} by enabling equilibria enumeration and selection. Compared to \gls{SGM}, \gls{eSGM} requires all variables to be integer-constrained (\ie, $k_i=0$), and it uses the mixed-integer optimization approach by Sandholm \etal~\cite{sandholm_mixed-integer_2005} in the play phase. \gls{eSGM} is only guaranteed to terminate if the strategy sets are bounded; this latter requirement ensures the existence of an equilibrium. In terms of equilibria selection, \gls{eSGM} selects the \emph{most probable} equilibrium based on the results of ~\cite{harsanyi_new_1995}.
    \item \gls{CNP}~\cite{Dragotto_2021_CNP} is an outer approximation algorithm to compute a Nash equilibrium for \emph{polyhedrally-representable separable games}, a large class of games contained in \glspl{IPG} where, for every player $i$, $f^i$ is a linear function of $x^i$, and $\cl \conv(\X{i})$ is a polyhedron. The authors generalize a result from~\cite{WNMS1}, and prove that the mixed-strategy sets $\Delta^i$ coincide with $\cl \conv(\X{i})$ whenever the players' objective functions are separable. Compared to the other algorithms we present, \gls{CNP} can handle unbounded (or uncountable) strategy sets. In its approximation phase, the algorithm approximates $\cl \conv(\X{i})$ with increasingly-accurate polyhedral outer approximations. The resulting approximated game consists of $n$ parametrized linear problems, as each player optimizes a linear function in $x^i$ over a polyhedral strategy set. Because of this characteristic, in the play phase, the algorithm solves a complementarity problem that simultaneously enforces the optimality conditions of the players' optimization problems. In the improvement phase, the algorithm performs both a membership step and a stability step and employs the resulting information to refine the outer approximations or return an equilibrium. The algorithm does not support equilibria enumeration and only implements some heuristics for equilibria selection.
    \item \gls{BM}~\cite{sagratella_computing_2016} is an outer approximation algorithm to compute and enumerate pure Nash equilibria in \glspl{IPG} with convex integer strategy spaces (\ie, convex constraints plus integer requirements) and continuously differentiable and concave payoff function $f^i$ in $x$. The algorithm explores the search spaces via a branching method that fixes variables and prunes subproblems not leading to equilibria. In its approximation phase, the algorithm approximates the game as a continuous game by relaxing the integrality requirements. In its play phase, the algorithm solves complementarity or variational inequality problems. In the improvement phase, the algorithm performs the membership step by enforcing that the tentative equilibrium fulfills integrality requirements and the stability step by computing the best response. The algorithm supports equilibria enumeration, but it cannot select equilibria.
    \item \gls{ZEROR}~\cite{ZERORegets} is an outer approximation algorithm to compute, select and enumerate Nash equilibria in \glspl{IPG} with mixed-integer linear strategy spaces and linearizable payoff functions, \ie, nonlinear payoff functions that can be reformulated as linear functions by adding extra variables and inequalities (\eg, see~\cite{vielma_mixed_2015}). This cutting-plane algorithm leverages the concept of equilibrium inequality, a linear inequality valid for any Nash equilibrium, and the associated improvement oracle. In its play phase, the algorithm solves a mixed-integer optimization problem that optimizes a given \emph{selection function} $h:\prod_{i=1}^n \X{i} \rightarrow \mathbb{R}$ over the joint strategy spaces of all players; naturally, the selection function should be computationally tractable in practice. In the improvement phase, the algorithm does not perform the membership step, as the result of the play phase is always an integer-feasible tentative equilibrium; however, it performs a stability step by computing the best responses and possibly generating equilibrium inequalities. The algorithm supports equilibria enumeration and selection and can also handle linear constraints that depend on $x$ (as opposed to $x^i$). In terms of equilibria selection, \gls{ZEROR} can select the equilibrium that maximizes the given selection function $h(x^1,\dots,x^n)$, as long as this function is computationally-tractable.
    \item \gls{BNP}~\cite{schwarze_branch-and-prune_2022} is an outer approximation algorithm to compute and enumerate pure Nash equilibria in \glspl{IPG} with convex integer strategy spaces and payoff functions $f^i$ that are strictly concave in the individual variables $x^i_j$ for any $j=1,\dots,k_i+n_i$. The algorithm consists of a branch-and-prune routine built on top of \gls{BM} with an additional pruning criterion. In the play phase, the algorithm only requires a feasible solution to the outer approximated optimization problems (\ie, it does not require a tentative equilibrium); however, the authors emphasize that computing a tentative equilibrium may improve the performance of their algorithm.
\end{enumerate}

\paragraph{Existence of Equilibria and Boundedness. } As highlighted in the analysis above, all the algorithms we review, except \gls{CNP}, require the players' strategy sets to be bounded. Indeed, boundedness is an instrumental assumption for the convergence of the algorithms to an equilibrium. If the strategy sets $\X{i}$ are bounded and integer-constrained, the \gls{IPG} has an equivalent normal-form game, and Nash's theorem on the existence of equilibria directly applies~\cite{nash_equilibrium_1950,nash_noncoop_1951}; this is the case for \gls{eSGM}, \gls{BM}, \gls{ZEROR} and \gls{BNP}. 
However, if the strategy sets $\X{i}$ include, as in \gls{CNP} and \gls{SGM}, continuous variables, the compactness of $\X{i}$ is required so that the existence result by Glicksberg holds~\cite{glicksberg_further_1952}. Moreover, \gls{CNP} and \gls{SGM} require the players' payoff functions to be separable in order to guarantee that any mixed equilibrium has finite support ~\cite{stein_separable_2008}.
Finally, \gls{CNP} is the only algorithm that allows unbounded strategy sets and returns an equilibrium or a proof of its non-existence. This peculiarity stems from the structure of the approximations \gls{CNP} builds. Specifically, whenever the players' payoff functions are separable, finding a mixed equilibrium is equivalent to finding a pure equilibrium in a \emph{convexified} version of the original game where each player strategy set is $\cl \conv(\X{i})$ as opposed to $\X{i}$. In practice, if one can build an algorithmic procedure to describe $\cl \conv(\X{i})$, then \gls{CNP} guarantees termination even with unbounded and mixed-integer strategy sets.

\section{An Example: Protecting Critical Infrastructure}
\label{sec:criticalnode}
This section provides an example of how to employ \glspl{IPG} to model practical problems and derive managerial insights from their solutions. Our example regards the protection of critical infrastructure of a network against external malicious attacks, for instance,  the protection of the population from a human-induced epidemic or of a computer network from a cyberattack. Specifically, we focus on determining the best strategies for the network operator (\eg, the public health officer or the network administrator) to protect their network from malicious attacks.

We compare two models: the \gls{CNG}, an \gls{IPG} inspired by \cite{Dragotto_critical_2023}, and the \gls{MCNP}, a bilevel model inspired by the sequential critical node problem of~\cite{baggio2021multilevel}. We remark that \gls{CNG} is a simultaneous game, whereas \gls{MCNP} is a sequential game, \ie, a game where players decide in a given order. Both models extend the critical node problem, \ie, the problem of detecting the most critical nodes in a graph for a given connectivity metric (see, \eg, \cite{lalou_critical_2018} for a survey), to the case where nodes can be protected and attacked. Our two game models differ fundamentally in terms of the interacting dynamics among the players. In~\cite{baggio2021multilevel}, the authors propose a sequential (\ie, trilevel) Defender-Attacker-Defender game, where the defender and attacker sequentially decide which nodes to protect and attack, respectively. On the contrary, \gls{CNG} assumes that the defender and the attacker simultaneously decide their strategies.

\subsection{Model Basics}
In the following, we describe the players, their strategies, and payoffs as defined in~\cite{Dragotto_critical_2023}. We then present the simultaneous and sequential dynamics and formally define the two mathematical models.

Two decision-makers, the attacker and the defender, decide how to protect a network represented by a finite set of resources $R$; these resources can, for instance, derive from the nodes of a graph representing the network topology. For any resource $i \in R$, we introduce a binary variable $x_i$ equal to $1$ if and only if the defender \emph{protects} resource $i$. Similarly, we introduce a binary variable $\alpha_i$ equal to $1$ if and only if the attacker \emph{attacks} resource $i$. The defender maximizes a function $\defender$ of $x$ parameterized in $\alpha$ representing the resources' availability (or the operativeness) under the strategies $x$ and $\alpha$. Symmetrically, the attacker maximizes a function $\attacker$ of $\alpha$ parameterized in $x$ representing the disruption on $R$ under $x$ and $\alpha$. The weights $a\in \mathbb{R}^{|R|}_+$ and $d\in \mathbb{R}^{|R|}_+$ represent the cost the player pays for selecting (\ie, defending or attacking) the resources, while $A \in \mathbb{R}_+$ and $D\in \mathbb{R}_+$ represent the players' overall budgets. Let $\pprotect_i$ and $\pattack_i$ represent the criticality of resource $i$ for the defender and the attacker, respectively. These parameters can quantify, for instance, the vulnerabilities the attacker could exploit or the importance of $i$ in the defender's network. We describe the functions $\defender$ and $\attacker$ in terms of $\pprotect_i$ and $\pattack_i$ according to whether resource $i$ is protected ($x_i=1$) or attacked ($\alpha_i=1$). Let $\delta$, $\eta$, $\epsilon$ and $\gamma$ be real-valued scalar parameters in $[0,1]$, so that $\delta < \eta < \epsilon$. The payoff contributions for each resource $i \in R$ are as follows:

\begin{enumerate}
    \item \textbf{Normal operations.} If $x_i=0$ and $\alpha_i=0$, the defender receives a payoff of $\pprotect_i$ as no attack is ongoing on $i$. The attacker pays an opportunity cost $\gamma \pattack_i$.
    \item \textbf{Successful attack.} If $x_i=0$ and $\alpha_i=1$, the attacker gets a payoff of $\pattack_i$, as the attacker successfully attacks resource $i$. The defender's operations on resource $i$ are impacted and worsen from $\pprotect_i$ to $\delta \pprotect_i$.
    \item \textbf{Mitigated attack.} If $x_i=1$ and $\alpha_i=1$, the defender mitigates the attempted attack on resource $i$. The defender's operations worsen from $\pprotect_i$ to $\eta \pprotect_i$, and the attacker gets a payoff of $(1-\eta) \pattack_i$.
    \item \textbf{Mitigation without attack.} If $x_i=1$ and $\alpha_i=0$, the defender protects resource $i$ without the attacker attacking it. The defender's operations worsen from $\pprotect_i$ to $\epsilon \pprotect_i$, and the attacker receives a payoff of $0$.
\end{enumerate}
All considered, the defender's payoff is
\begin{align}
    \defender(x,\alpha)=\sum_{i \in R} \Big (\pprotect_i\big ( (1-x_i)(1-\alpha_i) + \eta x_i\alpha_i  +\epsilon x_i(1-\alpha_i) + \delta(1-x_i)\alpha_i \big) \Big ),
\end{align}
and the attacker's payoff is
\begin{align}
    \attacker(\alpha,x)=\sum_{i \in R} \Big (\pattack_i\big ( -\gamma(1-x_i)(1-\alpha_i) + (1-x_i)\alpha_i  +(1-\eta)x_i\alpha_i \big) \Big ).
\end{align}

\paragraph{Simultaneous Game.} In the simultaneous \gls{CNG} of~\cite{Dragotto_critical_2023}, the attacker and the defender \emph{simultaneously} decide their strategies as in \cref{def:CNG}.

\begin{definition}[Critical Node Game~\cite{Dragotto_critical_2023}]
    The \gls{CNG} is a 2-player \emph{simultaneous} and \emph{non-cooperative} \gls{IPG} with \emph{complete information} where the first player (the defender) solves
    \begin{align}
        \underset{x}{\maxi} \Big \{ \defender(x;\alpha) : d^\top x \le D, x \in \{0,1\}^{|R|} \Big\},
    \end{align}
    and the second player (the attacker) solves
    \begin{align}
        \underset{\alpha}{\maxi} \Big \{ \attacker(\alpha;x) : a^\top \alpha \le A, \alpha \in \{0,1\}^{|R|} \Big\}.
    \end{align}
    \label{def:CNG}
\end{definition}
In the \gls{CNG}, the strategies of the attack and defense are uncertain, as the players decide simultaneously and are unaware of their opponent's strategy; precisely, the parameters associated with the attacker's optimization problem can model the \emph{class} of attacks for which the defender would like to develop protective strategies. In this sense, the simultaneousness of the \gls{CNG} enables the modeling of a wide range of attacks without resorting to any assumptions on the order of play for the players.

\vspace{4pt}

\paragraph{Sequential Game.} Symmetrically, in the \gls{MCNP} inspired from the bilevel variant of ~\cite[Section 3.1]{baggio2021multilevel}, the defender and the attacker \emph{sequentially} decide their strategies as in \cref{def:MCNP}.

\begin{definition}[Multilevel Critical Node Problem~\cite{baggio2021multilevel}]
    The \gls{MCNP} is a 2-player \emph{sequential} game formulated as the following bilevel program:
    \begin{align}
        \underset{\substack{x, \hat{\alpha}}}{\maxi} \quad &  \defender(x,\hat{\alpha}) \\
         \text{s.t} \quad & d^\top x \le D,\\
        & x \in \{0,1\}^{|R|}, \\        
        & \hat{\alpha} \in  \arg \max_{\alpha} \Big \{ \attacker(\alpha;x) : a^\top \alpha \le A, \alpha \in \{0,1\}^{|R|} \Big \}.
    \end{align}
    \label{def:MCNP}
\end{definition}
In contrast to the \gls{CNG}, in the \gls{MCNP}, the defender commits to a defensive strategy, and then, the attacker observes it and reacts by playing its optimal attack. The defender's and the attacker's constraints and objectives are analogous to the ones of the  \gls{CNG}.

\subsection{Experiments}
We solve the simultaneous \gls{CNG} of \cref{def:CNG} through the open-source solver ZERO~\cite{Dragotto_2021_ZEROSoftware} with \gls{ZEROR} and \gls{CNP}. The solver, the instances, and our implementations are available at \url{https://github.com/ds4dm/ZERO}. For this tutorial, we also provide an open-source implementation of  \gls{ZEROR} inside ZERO.
We implemented the algorithm by formulating the three building blocks (\ie, play, approximate, and improve) using the components available in ZERO. 

We run our experiments on $8$ CPUs and $64$ GB of RAM, with \emph{Gurobi 10} as the optimization solver. \footnote{Specifically, we employ $8$ Intel Xeon Gold $6142$ CPUs running at a speed of 2.60GHz each.} We generate a series of synthetic instances according to the generation scheme proposed in~\cite{Dragotto_critical_2023}. The instance generation attempts to obtain instances that mimic the dynamics of different types of realistic attacks in cloud networks; we refer to \cite{Dragotto_critical_2023} for a detailed description of the instances. In terms of cardinality, we consider problems with $10,20,25,50$ resources. For each cardinality of set $R$, we generate $20$ instances with different parameters.
We perform two series of experiments. First, we compare the performance of \gls{CNP} and \gls{ZEROR} for the problem of computing equilibria. Second, we qualitatively discuss the insights we can obtain by comparing the solutions of the \gls{CNG} to the ones of the \gls{MCNP}. 

We remark that neither of these two sets of experiments is intended to be exhaustive. In the first series, we aim to showcase that two different algorithms can solve \glspl{IPG}; the algorithms naturally have different characteristics, as discussed in this tutorial, and focus on different concepts of Nash equilibria. Furthermore, we aim to facilitate the development of new, improved algorithms for \glspl{IPG} by showcasing the use of the open-source code in ZERO~\cite{Dragotto_2021_ZEROSoftware} for our implementations. The second series of experiments aims to provide some (preliminary) insights into the different modeling paradigms for the multi-agent critical node problem (\eg, sequential or simultaneous) and the qualitative repercussions of the modeling paradigm on the solutions.

\vspace{4pt}

\paragraph{Performance Comparison.} First, we compare the performance of \gls{CNP} and \gls{ZEROR} for the task of computing equilibria.
We remark that \gls{CNP} can compute mixed equilibria, whereas \gls{ZEROR} can only compute pure equilibria. 
In this set of experiments, \gls{CNP} computes \emph{an} equilibrium (\ie, either a pure or mixed Nash equilibrium), whereas \gls{ZEROR} computes the pure equilibrium that maximizes the defender's payoff $\defender$; In practice, this implies that \gls{ZEROR} may determine that no pure equilibrium exists and return a certificate of non-existence or provide a pure equilibrium maximizing $\defender$; on the contrary, \gls{CNP} always returns an equilibrium, as there always exists a mixed equilibrium in the \gls{CNG} (since each player has finitely many strategies). 
We report a summary of our computational experiments in \cref{tab:performance}. Specifically, for each cardinality of $R$ and algorithm, we report: 
\begin{enumerate*}
\item the number of instances for which the algorithm computed one equilibrium, in particular, one pure equilibrium (\textit{\#Eq} and \textit{\#Pure eq}), 
\item the number of times the algorithm reached the time limit of $50$ seconds (\textit{\#TL}),
\item the average computing time in seconds (\textit{Time (s)}), and
\item the average number of iterations (\textit{\#Iterations}), \ie, runs of the three-phase process.
\end{enumerate*}

\begin{table}[!ht]
    \centering
    \caption{Performance comparison of \gls{CNP} and \gls{ZEROR}.\label{tab:performance}}
    \resizebox{0.8\textwidth}{!}{
        \begin{tabular}{c@{\hspace{1.5em}}l@{\hspace{3em}}r@{\hspace{1em}}r@{\hspace{1em}}r@{\hspace{1em}}r@{\hspace{1em}}r@{\hspace{1em}}}
            \hline
            $\bm{|R|}$ & \textbf{Algorithm} & \textbf{\#Eq} & \textbf{\#Pure eq} & \textbf{\#TL} & \textbf{Time (s)} & \textbf{\#Iterations} \\
            \hline
            10         & \gls{CNP}         & 20            & 16                 & 0          & 0.03              & 10                    \\
                     & \gls{ZEROR}        & 16            & 16                 & 0          & 0.03              & 3                     \\ \hline 
            20         & \gls{CNP}         & 20            & 8                  & 0          & 0.17              & 22                    \\
                     & \gls{ZEROR}        & 8             & 8                  & 2          & 5.41              & 14                    \\ \hline
            25         & \gls{CNP}         & 20            & 4                  & 0          & 0.55              & 50                    \\
                     & \gls{ZEROR}        & 4             & 4                  & 2          & 6.32              & 15                   \\ \hline
            50         & \gls{CNP}         & 20            & 0                  & 0          & 5.17              & 132                   \\
                     & \gls{ZEROR}        & 0             & 0                  & 8          & 21.46             & 27                    \\
            \hline
        \end{tabular}
    }
\end{table}

As expected, we observe that \gls{CNP} always computes an equilibrium within reasonable computing times and never reaches the time limit. On the contrary, \gls{ZEROR} tends to be slower than \gls{CNP}, reaching the time limit in $12$ instances. For a significant number of instances, \gls{ZEROR} certifies that no pure equilibrium exists, namely $4$, $10$, $14$, and $12$ for $|R| \in \{10, 20, 25, 50\}$, respectively. 
Interestingly, whenever \gls{ZEROR} returns a pure equilibrium, also \gls{CNP} returns a pure equilibrium.
Regarding the number of iterations, that of \gls{CNP} tends to be significantly higher than the iterations required by \gls{ZEROR} to find an equilibrium.
We speculate this mainly depends on the different types of equilibria the algorithms can compute. Whereas \gls{ZEROR} only computes pure equilibria, \gls{CNP} computes mixed Nash equilibria, thus considering a larger search space. We also remark that the algorithms have different play phases. On the one hand, \gls{CNP} solves a complementarity problem (\eg, with PATH \cite{dirkse_path_1995,ferris_interfaces_1999} whenever the complementarity problem is linear) and only requires a feasible solution. On the other hand, \gls{ZEROR} solves to optimality an integer program (with Gurobi). Empirically, the latter often requires a longer solution time compared to the former.
Finally, as the number of strategies increases, \gls{CNG} often does not have pure Nash equilibria, \eg, when $|R|=50$, suggesting that the randomization over a small set of strategies is more likely part of a rational strategy for the players.

\vspace{4pt}

\paragraph{Sequential vs. Simultaneous.} We perform a second series of experiments by qualitatively comparing the simultaneous \gls{CNG} with the sequential \gls{MCNP} defender-attacker bilevel problem to gain insights on their solutions. On the one hand, we solve the \gls{CNG} with \gls{ZEROR} and compute the (absolute approximate) pure equilibrium maximizing the defender's payoff, \ie, maximizing $\defender(x,\alpha)$; recall the definition of approximate equilibrium from Section~\ref{sec:what}. On the other hand, we solve the \gls{MCNP} with the bilevel solver of~\cite{fischetti_new_2017} by maximizing the defender's payoff (\ie, the leader's objective). For our tests, we employ a metric inspired from the \gls{POS} introduced in~\cite{Dragotto_critical_2023}, \ie, the ratio between the best defender's objective under any outcome of the game and the defender's objective under the equilibrium maximizing it. Formally, we define the \gls{POS} of a solution $(\bar x, \bar \alpha)$ as 
\begin{align}
\text{PoS}(\bar x, \bar \alpha) \quad = \quad \frac{\underset{x,\alpha}{\maxi} \{ \defender(x,\alpha) : d^\top x \le D, a^\top \alpha \le A, x \in \{0,1\}^{|R|}, \alpha \in \{0,1\}^{|R|}\}}{\defender(\bar x,\bar \alpha) },
\end{align}
where $\defender(\bar x,\bar \alpha)$ denotes the value (\ie, defender's payoff) of $(\bar x, \bar \alpha)$ of either the \gls{CNG} solved with \gls{ZEROR} or the \gls{MCNP} solved with the bilevel solver of~\cite{fischetti_new_2017}.
In practice, a solution with a low \gls{POS} intuitively indicates that the defensive strategy is more efficient, as the payoff of $(\bar x, \bar \alpha)$ is closer to those of the best strategies the defender would play if it could force the attacker to play any given strategy.
We report a summary of our computational experiments in \cref{tab:models}. Specifically, for each cardinality of $R$ and algorithm, we report: 
\begin{enumerate*}
\item the absolute approximation constant for the equilibrium (\textit{Epx}), 
\item the number of times the algorithm reached the time limit of $120$ seconds (\textit{\#TL}),
\item the average computing time in seconds (\textit{Time (s)}), and
\item the average \gls{POS} (\emph{PoS}).
\end{enumerate*}

\begin{table}[ht]
    \centering
    \caption{Comparison for the simultaneous \gls{CNG} (\gls{ZEROR}) and the sequential \gls{MCNP} (bilevel solver from~\cite{fischetti_new_2017}) .\label{tab:models}}
    \resizebox{0.75\textwidth}{!}{
        \begin{tabular}{c@{\hspace{1.5em}}l@{\hspace{3em}}l@{\hspace{3em}}r@{\hspace{1.5em}}r@{\hspace{1.5em}}r@{\hspace{1.5em}}r@{\hspace{1.5em}}}
            \hline
            $\bm{|R|}$ & \textbf{Model} &\textbf{Algorithm} & \textbf{Epx} & \textbf{\#TL} & \textbf{Time (s)} & \textbf{PoS} \\
            \hline
            10         & \gls{MCNP} & Bilevel             & -         & 0         & 0.07              & 1.43         \\
                       & \gls{CNG}  & \gls{ZEROR}        & 1.50         & 0         & 0.10              & 1.05         \\ \hline
            20         & \gls{MCNP} & Bilevel            & -         & 0         & 15.07             & 1.53         \\
                       & \gls{CNG}  & \gls{ZEROR}        & 4.40         & 2         & 20.27             & 1.12         \\ \hline
            25         & \gls{MCNP} & Bilevel            & -         & 20        & 120.00            & 1.36         \\
                       & \gls{CNG}  & \gls{ZEROR}        & 6.47         & 4         & 36.42             & 1.26          \\ \hline
            50         & \gls{MCNP} & Bilevel            & -         & 20        & 120.00            & 1.37         \\
                       & \gls{CNG}  & \gls{ZEROR}        & 17.50        & 8         & 72.56             & 1.14\\
            \hline
        \end{tabular}
    }
\end{table}
We observe that, in these instances, the \gls{POS} achieved by the \gls{CNG} (solved by \gls{ZEROR}) is consistently lower than that of the \gls{MCNP} (solved by the bilevel solver). This is due to the different modeling paradigms of the two models. The \gls{CNG} is indeed closer to the idea of robust optimization: instead of considering a single attack, the defender models the \emph{class} of attacks that could happen, given some minimal knowledge about the attacker's capabilities (\eg, the budget constraint and the importance of each resource). In contrast, the \gls{MCNP} models a sequential decision-making process in which the defender considers the attacker as a follower that adapts to the defender's decisions; this implies that the model assumes the attacker can witness the defender's decisions before committing to an attack strategy. In this sense, the \gls{CNG} provides a more realistic model for protecting critical resources under limited observability of the attacker's and defender's strategies.

Finally, we remark that we can always recover feasible solutions and equilibria even when the algorithms hit the time limit without proving optimality. Specifically, \gls{ZEROR} returns the equilibrium with the lowest absolute approximation constant, \ie, the equilibrium with the smallest \emph{Epx} found. Similarly, the bilevel solver returns the best incumbent solution found in the search process; we remark that finding a feasible solution is trivial in the latter case.

\section{Perspectives and Conclusions}
\label{sec:ideas}
We conclude our tutorial by presenting ideas regarding future research directions in \glspl{IPG}. Specifically, we will introduce the concept of generalized Nash equilibrium, discuss the assumption of complete information and parameter estimation, and present some algorithmic ideas. Finally, we conclude in \cref{sec:conclusions}. 

\subsection{Generalized Nash Equilibria Problems} 
In some practical context, the players' strategy sets $\X{i}$ may depend on the opponents' variables, \ie, $\X{i}(x^{-i})$ is parametrized in $x^{-i}$. For instance, consider the game in \cref{ex:generalized}.

\begin{example}[A Generalized \gls{IPG}]
Consider an \gls{IPG} where 
\begin{enumerate*}
\item each player needs to route an integral amount of traffic through a network defined over a graph,
\item each player incurs in a non-negative cost when selecting a specific edge, and 
\item each edge can be selected by at most one player
\end{enumerate*}. In this example, each player's strategy impacts its opponents' available strategies (\ie, the available edges). 
\label{ex:generalized}
\end{example}

Nash equilibria problems are simultaneous games where the strategy set of each player can depend on the opponents' strategies. Therefore, we call the extension of \glspl{IPG} to this form as \emph{generalized} \glspl{IPG}. The concept of Nash equilibrium of \cref{def:NE} still holds; however, the quantifier of \cref{eq:NE} becomes $\forall \tilde{x}^i \in \X{i}(x^{-i})$. Similarly, the equilibria, in this case, can be pure and approximate. We remark that, in this tutorial, we specifically focused on \glspl{IPG}, as the computational literature on generalized \glspl{IPG} is still not mature enough. However, from a practical perspective, \emph{generalized} \glspl{IPG} and their equilibria can provide a compelling modeling extension. In this context, we refer the reader to the works of Sagratella~\cite{sagratella_2019} and Harks and Schwarz~\cite{harks_generalized_2022}, where the authors propose two algorithms to compute equilibria in \emph{generalized} \glspl{IPG}.

\subsection{Incomplete Information and Parameters Estimation} 
One of the most restrictive assumptions in \glspl{IPG} is the one of complete information, namely, the assumption that the players know their optimization problems and those of their opponents. 
From an economic standpoint, players can, on the contrary, have different levels of information, and this should be reflected in the game; for instance, a player may have no information about some of the parameters describing the objective function or feasible set of its opponents.
We believe that the imperfect assumption of complete information is mainly due to the novelty of \glspl{IPG} and to the algorithmic challenges associated with the computation of their Nash equilibria. Indeed, with no exception, all the algorithms we review in this tutorial have been developed in the past decade. Further work is needed to drop the assumption of complete information and embrace incomplete information, \eg, by incorporating ideas and concepts from Bayesian games.
In \cref{sub:why}, we explain why \glspl{IPG} and their Nash equilibria are powerful design and intervention tools for regulators. However, from a practical perspective, assuming that regulators can directly observe the \emph{parameters} of the game (\ie, the payoff functions and strategy sets of the players) might be rather impractical in some settings.
Several works at the interface of machine learning and game theory recently addressed this issue by providing efficient methods to estimate these parameters starting from data containing observed Nash equilibria \cite{bertsimas_data-driven_2015,Dragotto_2023_LearningPotential,ling_what_2018,fabiani_learning-based_2022}. In this sense, the algorithms for computing equilibria in \glspl{IPG} are the fundamental tools to develop more sophisticated approaches where the parameters of the game are dynamically estimated over time.

\vspace{4pt}

\paragraph{Incomplete Information.} Through the estimation of parameters, we can also define \emph{incomplete information} \glspl{IPG} by adapting the well-known model of beliefs developed by Harsanyi~\cite{Harsanyi1987}. In this setting, an \gls{IPG} with incomplete information is an \gls{IPG} complemented with a set of \emph{types} for each player and a joint probability distribution over them; specifically, a type for a player represents a common belief on the payoff function for that player. We can also interpret this \gls{IPG} class using the language of stochastic optimization: the types correspond to \emph{scenarios}, \ie, a specification of a type for each player, and its probability of occurrence is a common belief of the players. Finally, this research stream intrinsically intersects with the concept of stochastic games, introduced by Shapely ~\cite{shapley1953stochastic}. In a stochastic game, the outcome of each move is determined by a random process that is not in the control of the players. We believe that exploring the interplay of mixed-integer stochastic optimization, stochastic games, and \emph{incomplete information} \glspl{IPG} is a natural direction for future work.

\subsection{Algorithmic Development}
We conclude this section by discussing ideas for future algorithmic developments in \glspl{IPG}. As mentioned, since the algorithms reviewed in this tutorial have been developed over the past decade, there is significant room for improvement, both from methodological and applied perspectives. 

\paragraph{Selecting Equilibria.} Only a few of the surveyed algorithms can optimize (or select) over the set of Nash equilibria, and further research in this direction is needed. Indeed, the selection is naturally a desirable property for an algorithm, especially when the \gls{IPG} instance has multiple equilibria. From a practical perspective, selecting equilibria enables regulators to privilege the equilibria exhibiting the most desirable properties and, on the other hand, apply equilibrium selection or refinement theories (see, \eg, ~\cite{harsanyi_general_1992}) to model realistic equilibrium choices by the players. 

\paragraph{Finding Better Approximations.} Most of the algorithms we presented require an approximated strategy space for each player in the first (approximation) phase. We speculate that better initial approximation strategies may positively impact the computational performance of the algorithms. This argument follows the logic of linear relaxations and integer programs: “tighter” (\ie, more accurate) starting linear relaxations generally reduce the computing times needed to solve the original integer program. Naturally, in the \gls{IPG} context, the concept of “tightness” is less clearly defined. In general, we believe that the initial approximations should include, as much as possible, the strategies the players are most likely to play in an equilibrium. In other words, the initial approximations should include the strategies that the players will \emph{rationally} play and exclude the others. A counter-intuitive aspect of this rationale is that some strategies that are optimal solutions to a player's optimization problem may never be played. This is the case, for instance, for strategies that are never best responses to any of the opponents' strategies.

\paragraph{Other Solution Concepts.} In the play phase, most algorithms compute tentative equilibria to the approximated game. In contrast, we suggest the more unconventional idea of computing different notions of equilibria, \eg, approximate equilibria, when such an approximation is faster to compute. This could provide a significant speedup in the play phase and help adapt the existing algorithms to different solution concepts other than the Nash equilibrium.

\paragraph{Releasing Open-source Software.} For the purpose of this tutorial, we implemented \gls{ZEROR} by formulating its building blocks in the open-source solver ZERO~\cite{Dragotto_2021_ZEROSoftware}. We believe that releasing open-source implementations of \gls{IPG} algorithms to the public can boost the interaction among researchers and further facilitate the development of new algorithms and impactful applications of \glspl{IPG}.

\subsection{Concluding Remarks}
\label{sec:conclusions}
In this tutorial, we presented a broad computational overview on computing Nash equilibria in \glspl{IPG}. This new class of games builds on top of the modeling capabilities of mixed-integer programming and enables decision-makers to account for their interaction with other decision-makers. We argued that \glspl{IPG} are powerful tools for their modeling capabilities, regulatory applications, and practical use cases. We critically reviewed the algorithmic literature on \glspl{IPG} and proposed a new taxonomy of the main components employed by the available algorithms. Furthermore, we provided an example in the context of protecting critical infrastructure and compared the informative content of the Nash equilibria with the solution of an analogous sequential model. Finally, we proposed a few ideas for future research directions, hoping to stimulate new interest from the operations research community.

\section*{Acknowledgments}
We are indebted to two anonymous reviewers whose detailed and insightful comments significantly helped to improve the paper.

\iftoggle{ARXIV}{%
\bibliographystyle{gabri}
}{
\bibliographystyle{template/TutORials}
}
\bibliography{biblio}

\end{document}